# On the Density of Rational and Integral Points on Algebraic Varieties

*Per Salberger*

## 0. Introduction

This paper is concerned with the number $N(X;B)$ of rational points of height at most $B$ on projective varieties $X \subset \mathbf{P}^n$ over $\mathbf{Q}$. To define the height $H(x)$ of a rational point $x$ on $\mathbf{P}^n$, we choose a primitive integral $n+1$-tuple $(x_0, ..., x_n)$ representing $x$ and let $H(x) = \max(|x_0|, ..., |x_n|)$. Our main result is the following.

**Theorem 0.1.** *Let $X \subset \mathbf{P}^n$ be a geometrically integral projective variety over $\mathbf{Q}$ of dimension $r$ and degree $d \geq 4$. Suppose that there are only finitely many $(r-1)$-planes over $\overline{\mathbf{Q}}$ on $X$. Then*

$$N(X;B) = O_{d,n,\varepsilon}(B^{r+\varepsilon}).$$

More generally it has been conjectured by Heath-Brown [BwHeSa], conj. 2 that the same uniform bound should hold for any geometrically integral projective variety of degree $d \geq 2$. The conjecture is established for $d \geq 6$ in [BwHeSa], cor.2 and it is also known when $d=2$ by [He$_2$], th. 2. One cannot expect a lower exponent in general as $N(X;B) \gg_X B^r$ whenever there is an $(r-1)$-plane on $X$ defined over $\mathbf{Q}$.

On using proposition 8.3, we obtain the following corollary of theorem 0.1.

**Corollary 0.2.** *Let $X \subset \mathbf{P}^{r+1}$ be a geometrically integral hypersurface of degree $d \geq 4$. Suppose that the singular locus of $X$ is of codimension at least three in $X$. Then*

$$N(X;B) = O_{d,r,\varepsilon}(B^{r+\varepsilon}).$$

The same assertion has been proved for non-singular hypersurfaces of degree $\geq 2$ in a series of papers of Heath-Brown and Browning (cf. [He$_1$],[He$_2$],[Bw], [BwHe$_2$], [BwHe$_3$]). To achieve this, they use exponential sums (mod $p$) when $r \geq 8$ and lattice reduction when $r \leq 7$. Our proof of corollary 0.2 is different and does not rely on these techniques.

The following result is a consequence of theorem 0.1 and proposition 9.3.

**Corollary 0.3.** *Let $X \subset \mathbf{P}^{r+2}$ be a geometrically integral r-dimensional intersection of two hypersurfaces of degree $d_1$ and $d_2$. Suppose that $d = d_1 d_2 \geq 4$ and that the singular locus of $X$ is of codimension at least three in $X$. Then*

$$N(X;B) = O_{d,r,\varepsilon}(B^{r+\varepsilon}).$$

This result seems to be new even for non-singular intersections of two quadrics.

To prove 0.1, we use the same approach as in [BwHeSa]. We first reduce to the case where $X$ is a hypersurface by means of a projection from $\mathbf{P}^n$ to $\mathbf{P}^{r+1}$ which maps $X$ birationally onto its image. Then $X \subset \mathbf{P}^{r+1}$ is defined by a form $F(x_0, ..., x_{r+1})$, which may be chosen such that $F$ is primitive with integer coefficients. Such a form defines also an affine scheme $Y \subset \mathbf{A}_{\mathbf{Z}}^{r+2}$. Now as any rational point on $X$ can be represented by an integral point on $Y$, we obtain theorem 0.1 for hypersurfaces from the following result.

**Theorem 0.4**. *Let $f(y_1,\ldots, y_n) \in \mathbf{Z}[y_1,\ldots, y_n]$, $n \geq 3$ be a polynomial of degree $d \geq 4$ which is absolutely irreducible over $\mathbf{Q}$ and let $n(Y; B)$ be the number of n-tuples $\mathbf{y} = (y_1,\ldots, y_n)$ of integers such that $y_1,\ldots, y_n \in [-B, B]$ and $f(\mathbf{y}) = 0$. Suppose that there are only finitely many affine linear $(n-2)$-subspaces over $\overline{\mathbf{Q}}$ on the affine hypersurface in $\mathbf{A}^n$ defined by $f(\mathbf{y}) = 0$. Then*
$$n(Y; B) = O_{d,n,\varepsilon}(B^{n-2+\varepsilon}).$$

Note that we allow $f$ to be a non-homogeneous polynomial. This makes it possible to follow the inductive approach with hyperplane sections in [Pi$_1$] and [BwHeSa]. But the induction step is much deeper than in (op.cit.) and related to the results of B.Segre [Sg] on varieties with many linear subspaces.

To begin the induction proof of theorem 0.4, we count integral points on affine surfaces by means of the determinant method of Heath-Brown [He$_2$], th.14. By using $p$-adic congruences between the coordinates of the integral points on $Y$, we construct auxiliary hypersurfaces of low degree such that any point of bounded height lies on one of these hypersurfaces. It is then enough to count points on the intersections of the original surface with the auxiliary hypersurfaces. To count points on the resulting curves, we construct a second set of auxiliary hypersurfaces by means of $q$-adic congruences for another prime $q$ and apply the theorem of Bézout.

To get a second set of auxiliary hypersurfaces which is as small as possible, we use the $p$-adic congruence conditions for the first prime $p$ once more. This is quite delicate since the points may specialise to singular $\mathbf{F}_p$-points on the curves. One can thus not use the implicit function theorem as in the proof of [He$_2$], th.14. To overcome this, we make a careful study of the stalks of the $\mathbf{F}_p$-points on models of the curves over $\mathbf{Z}$ and of their multiplicities on the reduction of the models (mod $p$). We use thereby classical theorems in algebraic geometry of M. Noether [Fu$_2$] and Kronecker-Castelnuovo [Ca].

This paper is organised as follows. In the first three sections we give a very general version of the determinant method in [He$_2$] and construct the auxiliary hypersurfaces mentioned above. In section 4 we use these auxiliary hypersurfaces to count points on space curves. In section 5 we collect some results on the geometry of algebraic surfaces in $\mathbf{P}^3$. These are used in section 6 in combination with the estimates for curves in sections 4 to prove theorem 0.4 for affine surfaces by means of the determinant method. In section 7 we collect some results on the geometry of varieties with many linear subspaces, which we then apply in section 8 in the proof of theorem 0.4 for affine hypersurfaces of dimension >2. Finally, in section 9 we deduce theorem 0.1 from theorem 0.4.

We shall in this paper use the convention that $\varepsilon$ may change between occurrences. We may thus, for example, first write $f(B)=O_\varepsilon(B^{2\varepsilon})$ and then $f(B)=O_\varepsilon(B^\varepsilon)$. Let us also point out that all intersections with projective linear subspaces are scheme-theoretical. The hyperplane section of an integral variety may thus be a reducible or non-reduced scheme.

# 1. Estimates of certain determinants

An important new tool in the study of counting functions of rational points of bounded height is the determinant method. It was initiated by Bombieri and Pila [BP] who used it to give uniform bounds for affine plane curves. Then, Heath-Brown gave a *p*-adic version of the method and developed it to a powerful tool for projective hypersurfaces (see [He$_2$], th.14). Finally, Broberg [Br] generalised the method in [He$_2$] to arbitrary projective varieties. The goal of the method in the projective case is to construct auxiliary forms of low degree containing all the rational points of bounded height that one would like to count. The existence of these auxiliary forms follows from the vanishing of certain determinants. We shall in this section estimate the absolute values of these determinants.

**Definition 1.1.** Let $K$ be a field and $I \subset K[x_0,..., x_n]$ be a homogeneous ideal. Then the Hilbert function $h_I(k) = \dim_K(K[x_0,...,x_n]/I)_k$. If $X \subset \mathbf{P}_K^n = \mathrm{Proj}(K[x_0,...,x_n])$ is a closed subscheme, then $h_X(k) = h_I(k)$ for the saturated homogeneous ideal $I$ corresponding to $X = \mathrm{Proj}(K[x_0,...,x_n]/I)$ (cf. p.125 in [Ha]).

It is well known (cf. pp. 51-52 in [Ha]) that there is a unique polynomial $p_I(t) \in \mathbf{Q}[t]$, the Hilbert polynomial, such that $h_I(k) = p_I(k)$ for all sufficiently large positive integers and that

(1.2) $\qquad p_I(t) = dt^r/r! +$ terms of lower degree ,

where $r$ is the dimension of $X = \mathrm{Proj}(K[x_0,...,x_n]/I)$ and $d$ is the degree of $X \subset \mathbf{P}^n$.

We shall also write $h(k)$ and $p(t)$ when it is clear which ideal $I$ we refer to.

The following two lemmas will be needed in order to prove that our estimates are uniform with respect to other parameters than $d$ and $n$.

**Lemma 1.3.** *Let $n$, $K$, $X = \mathrm{Proj}(K[x_0,...,x_n]/I)$ be as in* (1.1). *Then $X$ is defined by forms in $K[x_0,...,x_n]$ of degree bounded in terms of $n$ and $p_I(t)$.*

*Proof.* Let $H$ be the Hilbert scheme [G$_1$] and $W \subset H \times \mathbf{P}^n$ be a universal family of closed subschemes of $\mathbf{P}^n$ with Hilbert polynomial $p(t)$. There exists by (op.cit.) a closed immersion $H \subset \mathbf{P}^m$. The composite embedding $W \subset H \times \mathbf{P}^n \subset \mathbf{P}^m \times \mathbf{P}^n$ is defined by finitely many bihomogeneous polynomials. The assertion follows by specialising these bihomogeneous polynomials at the Hilbert point of $X$ on $H \subset \mathbf{P}^m$.

**Lemma 1.4.** *Let $d$ and $n$ be non-negative integers and $K$ be a field. Then there are only finitely many possibilities for the Hilbert function $h_X : \mathbf{Z}_+ \to \mathbf{Z}$ of closed geometrically reduced equidimensional subschemes $X \subset \mathbf{P}_K^n$ of degree $d$.*

*Proof.* See [Kl$_1$], cor. 6.11 or [RTV], th. 3.4 for a proof. (The result was first stated in [G$_1$], lemme 2.4 as pointed out to me by B. Conrad.)

We now follow the approach with Gröbner bases in Broberg's thesis [Br].

**Definition 1.5.** A graded monomial ordering $<$ is a total ordering on the set of monomials

$x^{\alpha} = x_0^{\alpha_0} \ldots x_n^{\alpha_n}$ satisfying the following conditions for the corresponding ordering on the set of exponents $\alpha = (\alpha_0, \ldots, \alpha_n) \in \mathbf{Z}_{\geq 0}^{n+1}$.

(i) $\alpha \geq \mathbf{0} = (0, \ldots, 0)$ for all $\alpha \in \mathbf{Z}_{\geq 0}^{n+1}$.

(ii) If $\alpha, \beta, \gamma \geq \mathbf{0} \in \mathbf{Z}_{\geq 0}^{n+1}$ and $\alpha+\gamma < \beta+\gamma$, then $\alpha < \beta$.

(iii) If $\alpha, \beta \in \mathbf{Z}_{\geq 0}^{n+1}$ and $\alpha \leq \beta$, then $\alpha_0 + \ldots + \alpha_n \leq \beta_0 + \ldots + \beta_n$.

The standard example is the *graded lexicographical ordering*. For this, $\alpha < \beta$ if and only if $\alpha_0 + \ldots + \alpha_n < \beta_0 + \ldots + \beta_n$ or $\alpha_0 + \ldots + \alpha_n = \beta_0 + \ldots + \beta_n$ and the left-most entry of $\alpha - \beta$ is negative.

**Definition 1.6.** The *leading monomial* of a polynomial $f(x_0, \ldots, x_n) = \sum a_{\alpha} x^{\alpha}$ with respect to a graded monomial ordering $<$ is the maximum of all $x^{\alpha}$ where $a_{\alpha} \neq 0$.

**Definition 1.7.** Let $K$ be a field and $I \subset K[x_0, \ldots, x_n]$ be a homogeneous ideal. Let $<$ be a graded monomial ordering on $\{x^{\alpha}, \alpha \in \mathbf{Z}_{\geq 0}^{n+1}\}$ and $m \in \{0, \ldots, n\}$. Then

$$\sigma_{I,m}(k) := \sum_{\alpha} \alpha_m$$

where $\alpha = (\alpha_0, \ldots, \alpha_n)$ runs over the exponent set of all monomials $x^{\alpha}$ which are not leading monomials of any form of degree $k$ in $I$. We shall also write $\sigma_m(k)$ instead of $\sigma_{I,m}(k)$ when it is clear which ideal $I$ we refer to.

The monomials of degree $k$ which are not leading monomials of any form of degree $k$ in $I$ form a $K$-basis of $(K[x_0, \ldots, x_n]/I)_k$ (see p. 452 in [CLO]). Hence

(1.8) $\quad\quad \sigma_{I,0}(k) + \ldots + \sigma_{I,n}(k) = k h_I(k)$.

**Lemma 1.9.** Let $K$, $I$, $<$, $\sigma_{I,m}$ be as above and suppose that $X = \mathrm{Proj}(K[x_0, \ldots, x_n]/I)$ is reduced and equidimensional. Then there exists a positive integer $k_0$ depending only on $n$ and $d = \deg X$ such that $\sigma_{I,m}(k)$, $m \in \{0, \ldots, n\}$ is given by a polynomial for $k \geq k_0$. In particular, there exists $a_{I,m} \geq 0$, $m \in \{0, \ldots, n\}$ such that

$$\sigma_{I,m}(k)/k h_I(k) = a_{I,m} + O_{d,n}(1/k), \quad\quad m \in \{0, \ldots, n\}.$$

*Proof.* It follows from lemma 1.3 and lemma 1.4 that $X$ is defined by forms of degree $<<_{d,n} 1$. The lemma is therefore a consequence of [Br], lemma 1.

**Remarks 1.10.** (a) $a_{I,0} + \ldots + a_{I,n} = 1$ by (1.8). Hence $0 \leq a_{I,m} \leq 1$ for $m \in \{0, \ldots, n\}$.

(b) If $I$ is generated by a form $F(x_0, \ldots, x_n)$ of degree $d$ with leading monomial
$x^{\alpha} = x_0^{\alpha_0} \ldots x_n^{\alpha_n}$, then $a_{I,m} = (d - \alpha_m)/nd$ for $m \in \{0, \ldots, n\}$.

The following result is a slight improvement of lemma 3 in [Br].

**Main lemma 1.11.** *Let $X \subset \mathbf{P}^n$ be a closed geometrically reduced equidimensional subscheme of dimension r and degree d and I be the saturated homogeneous ideal of $\mathbf{Q}[x_0,..., x_n]$ defined by X. Let $<$ be a graded monomial ordering on $\{x^\alpha, \alpha \in \mathbf{Z}_{\geq 0}^{n+1}\}$. Let $(B_0,..., B_n) \in \mathbf{R}_{\geq 1}^{n+1}$, and $\xi_l = (\xi_{l,0},..., \xi_{l,n})$, $l \in \{1, ..., s\}$ be (n+1)-tuples of integers representing rational points on X such that $|\xi_{l,m}| \leq B_m$ for all $l \in \{1, ..., s\}$ and all $m \in \{0, ..., n\}$. Then there exist monomials $F_1(x_0,..., x_n),..., F_s(x_0,..., x_n)$ of the same degree $k = (r!/d)^{1/r} s^{1/r} + O_{d,n}(1)$ which are not leading monomials of any form in I such that*

$$|\det(F_j(\xi_l))| \leq s!(B_0^{a_0} ... B_n^{a_n})^{kh(k)}(B_0 ... B_n)^\kappa,$$

*for some $\kappa = O_{d,n}(h(k))$.*

*Proof.* Let $k \geq 1$ be the unique integer such that $h(k-1) < s \leq h(k)$. Then $s = dk^r/r! + O_{d,n}(k^{r-1})$ and $k = (r!/d)^{1/r} s^{1/r} + O_{d,n}(1)$ as $h(k) = dk^r/r! + O_{d,n}(k^{r-1})$ by (1.2) and lemma 1.4.

To estimate the determinant, we use the trivial estimate

$$|\det(F_j(\xi_l))| \leq s! B_0^{\sigma_0(k)} ... B_n^{\sigma_n(k)}$$

and the fact (see lemma 1.9) that $\sigma_m(k) = a_m k h(k) + O_{d,n}(h(k))$ for $m \in \{0, ..., n\}$.

In order to apply this result to affine varieties, we shall need the following lemma.

**Lemma 1.12.** *Let $X \subset \mathbf{P}^n$ be a closed equidimensional subscheme of dimension r over a field K and $I \subset K[x_0,..., x_n]$ be a homogeneous ideal defining X. Suppose that the hyperplane $\Pi_0 \subset \mathbf{P}^n$ defined by $x_0 = 0$ intersects X properly. Then there exists a graded monomial ordering on the monomials in $x_0,..., x_n$ such that $a_{I,1} + ... + a_{I,n} \leq r/(r+1)$. In particular, if $(B_0,..., B_n) = (1, B,..., B)$, then $B_0^{a_0} ... B_n^{a_n} \leq B^{r/(r+1)}$.*

*Proof.* We shall use the *reverse* graded lexicographical ordering, where $\alpha < \beta$ if and only if $\alpha_0 + ... + \alpha_n < \beta_0 + ... + \beta_n$ or $\alpha_0 + ... + \alpha_n = \beta_0 + ... + \beta_n$ and the left-most entry of $\alpha - \beta$ is positive. Let $M_1,..., M_{h(k)}$ be the set of all monomials in $(x_0,..., x_n)$ of degree k which are not leading monomials of any form of degree k in I and $m_j(x_1,..., x_n) = M_j(1, x_1,..., x_n)$ for $j=1,..., h(k)$. Then, $\sigma_{I,1}(k) + ... + \sigma_{I,n}(k) := \deg m_1 + ... + \deg m_{h(k)}$. We must therefore show that :

$$\limsup_{k \to \infty} (\deg m_1 + ... + \deg m_{h(k)})/kh_I(k) \leq r/(r+1).$$

As $\lim_{k \to \infty} h_I(k)/dk^r = 1/r!$, $d = \deg X$, it is enough to show that

$$\deg m_1 + ... + \deg m_{h(k)} \leq dk^{r+1}/(r-1)!(r+1) + O_I(k^r)$$

when $k \to \infty$.

Let $J \subset K[x_0,..., x_n]$ be the homogeneous ideal generated by $x_0$ and I. Then, $m_j$ cannot be the leading monomial of a form F in J, since otherwise $M_j$ would be a leading monomial of $x_0^{k-\deg m_j} H \in I$ if $F = x_0 G + H$, $H \in I$. As the monomials $m_j$, $j=1,..., h_I(k)$ are different of degree at most k, we therefore get :

$$\deg m_1 + ... + \deg m_{h(k)} \leq h_J(1)1 + ... + h_J(k)k.$$

Also, $h_J(k)k \leq dk^r/(r-1)! + O_J(k^{r-1})$ by (1.2) since the scheme $\Pi_0 \cap X$ defined by $J$ is of the same degree as $X$. To complete the proof, use the fact that $1^r + \ldots + k^r = k^{r+1}/(r+1) + O_r(k^r)$.

## 2. *p*-power factors of certain determinants

We shall in this section (see 2.5) prove that the determinants in lemma 1.11 are divisible by high *p*-powers in case we consider integral points in the same congruence class (mod *p*). This will be used in the next section to show that the determinants vanish if *p* is large enough with respect to the maximum height of the rational points. Our results generalise the *p*-adic estimates in [He$_2$] and [Br] for non-singular congruence classes (mod *p*) to general congruence classes (mod *p*). We shall use the following notation throughout the section.

**Notation 2.1**. (a) $X \subset \mathbf{P}^n_\mathbf{Q} = \text{Proj}(\mathbf{Q}[x_0, \ldots, x_n])$ is a closed, reduced, equidimensional subscheme of dimension *r* and degree *d*.
(b) $I \subset \mathbf{Q}[x_0, \ldots, x_n]$ is the homogeneous ideal of polynomials which vanish at *X*.
(c) $\Xi = \text{Proj}(\mathbf{Z}[x_0, \ldots, x_n]/(I \cap \mathbf{Z}[x_0, \ldots, x_n]))$ is the scheme-theoretic closure of *X* in $\mathbf{P}^n_\mathbf{Z} = \text{Proj}(\mathbf{Z}[x_0, \ldots, x_n])$.
(d) $I_p \subset \mathbf{F}_p[x_0, \ldots, x_n]$ is the image of $I \cap \mathbf{Z}[x_0, \ldots, x_n]$ in $\mathbf{F}_p[x_0, \ldots, x_n]$.
(e) $X_p = \text{Proj}(\mathbf{F}_p[x_0, \ldots, x_n]/I_p) = \Xi \times_\mathbf{Z} \mathbf{F}_p$.

**Lemma 2.2**. *Let d and n be given and X, $X_p$ be as in 2.1(a). Then $X_p$ is defined by forms of degree $\ll_{d,n} 1$.*

*Proof.* Let *K* be an algebraic field extension of **Q**. Then $X_K$ is generically reduced and has no embedded components. Hence *X* is geometrically reduced so that we may apply lemma 1.4. There are thus $\ll_{d,n} 1$ possible Hilbert polynomials for the schemes $X \subset \mathbf{P}^n_\mathbf{Q}$ in 2.1(a). But the Hilbert polynomials of $X = \Xi \times_\mathbf{Z} \mathbf{Q}$ and $X_p = \Xi \times_\mathbf{Z} \mathbf{F}_p$ coincide since $\Xi \to \text{Spec } \mathbf{Z}$ is flat (see III.9.8 and III.9.9 in [Ha]). There are thus $\ll_{d,n} 1$ possible Hilbert polynomials also for the $\mathbf{F}_p$-subschemes $X_p \subset \mathbf{P}^n_{\mathbf{F}_p}$ in 2.1. The result now follows from lemma 1.3.

**Lemma 2.3**. *Let d and n be given and X, $X_p$ be as in 2.1(a). Let A be the stalk of $X_p$ at some $\mathbf{F}_p$-point P on $X_p$ of multiplicity $\mu$ and m be the maximal ideal of A. Let $g_{X,P} : \mathbf{Z}_+ \to \mathbf{Z}$ be the function where $g_{X,P}(k) = \dim_{A/m} m^k/m^{k+1}$. Then there are $\ll_{d,n} 1$ different functions $g_{X,P}$ among all pairs (X, P) as above. In particular, $g_{X,P}(k) = \mu k^{r-1}/(r-1)! + O_{d,n}(k^{r-2})$.*

*Proof.* $g_{X,P}$ is equal to the Hilbert function of the projectivised tangent cone $W_P \subset \mathbf{P}^{n-1}_{\mathbf{F}_p}$ of $X_p$ at *P* and $W_P \subset \mathbf{P}^{n-1}_{\mathbf{F}_p}$ is defined by forms of degree $\ll_{d,n} 1$ since this is true for $X_p$ (cf. III.3 in [Mu$_2$]). We may thus just as in the proof of [Br], lemma 1 choose a Gröbner base of forms of degree $\ll_{d,n} 1$ for the homogeneous ideal of $W_P \subset \mathbf{P}^{n-1}_{\mathbf{F}_p}$. The Hilbert function will not change (cf. p.452 in [CLO]) if we replace this set of forms by their leading monomials. There are thus $O_{d,n}(1)$ Hilbert functions for subschemes $W \subset \mathbf{P}^n$ defined by forms of degree $\ll_{d,n} 1$ and $O_{d,n}(1)$ Hilbert functions $g_{X,P}$ for the set of all pairs (X, P) as above. In particular, $g_{X,P}(k) = \mu k^{r-1}/(r-1)! + O_{d,n}(k^{r-2})$ by (1.2).

**Lemma 2.4**. *Let R be a commutative noetherian local ring containing $\mathbf{Z}_p$ as a subring and let $A=R/pR$. Let m be the maximal ideal of A and $(n_l(A))_{l=1}^{\infty}$ be the non-decreasing sequence of integers $k \geq 0$ where k occurs exactly $\dim_{A/m} m^k/m^{k+1}$ times. Let $r_1,...,r_s$ be elements of R and $\phi_1,...,\phi_s$ be ring homomorphisms from R to $\mathbf{Z}_p$. Then the determinant of the $s \times s$-matrix $(\phi_i(r_j))$ is divisible by $p^{A(s)}$ for $A(s) = n_1(A)+ n_2(A)+ ... + n_s(A)$.*

*Proof.* Let $I$ be the kernel of $\phi_1$. Then $I$ is sent to $m$ under $R \to R/pR=A$ and the induced ring homomorphism from $R/I$ to $A/m$ is just the map from $\mathbf{Z}_p$ to $\mathbf{F}_p$. The quotient $I^k/I^{k+1}$, $k \geq 0$ is thus a $\mathbf{Z}_p$-module generated by at most $\dim_{A/m} m^k/m^{k+1}$ elements by Nakayama's lemma.

Suppose there are more than $g(k) = \dim_{A/m} m^k/m^{k+1}$ elements in $I^k$-$I^{k+1}$ among $r_1,...,r_s$, say $r_1,...,r_q$. There is then a homomorphism of $\mathbf{Z}_p$-modules $\psi: \mathbf{Z}_p^q \to I^k/I^{k+1}$ where $\psi(\beta_1,...,\beta_q) = \beta_1 r_1+...+\beta_q r_q$ (mod $I^{k+1}$). But then as $I^k/I^{k+1}$ is generated by $g(k)<q$ elements we conclude that there is an element $(\beta_1,...,\beta_q)$ in the kernel of $\psi$ with one coordinate, say $\beta_q=1$. Then $\rho_q = \beta_1 r_1 + ... \beta_{q-1} r_{q-1} + r_q \in I^{k+1}$ and we will not change the determinant of $(\phi_i(r_j))$ if we replace $r_q$ by $\rho_q$. If we continue to make such elementary transformations, we will finally be in a situation where there are at most $\dim_{A/m} m^k/m^{k+1}$ elements in $I^k$-$I^{k+1}$ among $r_1,...,r_s$ for each $k \geq 0$. After a reordering we get that $r_j$ belongs to the $n_j(A)$-th power of $I$ for each $j=1,...,s$. This implies in its turn that all the elements $\phi_1(r_j),...,\phi_s(r_j)$ are divisible by the $n_j(A)$-th power of $p$ so that $\det(\phi_i(r_j))$ is divisible by $p^{A(s)}$.

**Main lemma 2.5**. *Let p be a prime and $X$, $r$, $d$, $n$, $\Xi$, $X_p$ be as in 2.1(a). Let P be an $\mathbf{F}_p$–point of multiplicity $\mu$ on $X_p$ and $\xi_1,...,\xi_s$ be primitive $(n+1)$-tuples of integers representing $\mathbf{Z}$-points on $\Xi$ with reduction P. Let $F_1,...,F_s$ be forms in $(x_0,...,x_n)$ with integer coefficients and $\det(F_j(\xi_l))$ be the determinant of the $s \times s$-matrix $(F_j(\xi_l))$. Then there exists a non-negative integer $N=(r!/\mu)^{1/r}(r/(r+1))s^{1+1/r} + O_{d,n}(s)$ such that $p^N | \det(F_j(\xi_l))$.*

*Proof.* We may after a coordinate change assume that the $x_0$-coordinates of $\xi_1,...,\xi_s$ are not divisible by $p$ and replace the forms $F_j$, $1 \leq j \leq s$ by the rational functions $f_j = F_j(1, x_1/x_0,..., x_n/x_0)$ without changing the $p$-adic valuation of the determinant. These rational functions are elements of $R$ and we get ring homomorphisms $\phi_1,...,\phi_s$ from $R$ to $\mathbf{Z}_p$ by evaluating at $\xi_1,...,\xi_s$. Hence $p^{A(s)} | \det(F_j(\xi_l))$ by the previous lemma. To complete the proof, we now show that $A(s) = (r!/\mu)^{1/r}(r/(r+1))s^{1+1/r} + O_{d,n}(s)$.

Let $g = g_{X,P} : \mathbf{Z}_{\geq 0} \to \mathbf{Z}$ be as in lemma 2.3 with $g(0)=1$ and $G(k)=g(0) + ... + g(k)$. Then from $g(k)= \mu k^{r-1}/(r-1)! + O_{d,n}(k^{r-2})$ we conclude that $G(k)= \mu k^r/r! + O_{d,n}(k^{r-1})$ and

$$(r!/\mu)^{1/r} (r/(r+1))G(k)^{1+1/r} = \mu k^{r+1}/(r-1)! (r+1) + O_{d,n}(k^r),$$

$$A(G(k)) = 0g(0) + ... + kg(k) = \mu k^{r+1}/(r-1)! (r+1) + O_{d,n}(k^r).$$

Hence $A(G(k)) = (r!/\mu)^{1/r}((r+1)/r)G(k)^{1+1/r} + O_{d,n}(G(k))$ as $k^r = O_{d,n}(G(k))$.

To deduce that $A(s) =(r!/\mu)^{1/r} (r/(r+1))s^{1+1/r} +O_{d,n}(s)$, let $k \geq 1$ be the unique integer such that $G(k-1)< s \leq G(k)$ and make use of the facts that

$$0 \leq A(G(k)) - A(s) \leq kg(k) <<_{d,n} k^r <<_{d,n} s,$$

$$0 \leq G(k)^{1+1/r} - s^{1+1/r} \leq G(k)^{1+1/r} - G(k-1)^{1+1/r} \ll_{d,n} k^r \ll_{d,n} s.$$

**Remark 2.6**. The case $\mu=1$ where $P$ is non-singular can be treated more directly by means of the implicit function theorem as in the proofs of [He$_2$], th.14 and [Br], lemma 3. But we shall also in the sequel need the more difficult case $\mu >1$ where $P$ is singular in order to improve on earlier estimates.

## 3. The auxiliary hypersurfaces

We shall in this section generalise the determinant method in [He$_2$] and [Br]. This produces auxiliary hypersurfaces for sets of rational points of bounded height on a projective variety satisfying congruence conditions. But we shall not as in (op.cit.) restrict to congruence classes corresponding to non-singular $\mathbf{F}_p$-points at one prime $p$.

The following notation will be used in the sequel.

**Notation 3.1**. Let $X \subset \mathbf{P}^n_{\mathbf{Q}} = \text{Proj}(\mathbf{Q}[x_0,..., x_n])$ be a closed subscheme over $\mathbf{Q}$ and let $\Xi$ be the scheme-theoretic closure of $X$ in $\mathbf{P}^n_{\mathbf{Z}} = \text{Proj}(\mathbf{Z}[x_0,..., x_n])$. Let $p_1,..., p_t$ be primes and $P_i$ be an $\mathbf{F}_{p_i}$-point on $X_{p_i} = \Xi \times_{\mathbf{Z}} \mathbf{F}_{p_i}$ for each $i \in \{1,..., t\}$. Let $U$ be a locally closed subvariety of $X$ defined over $\mathbf{Q}$. Then,

(i) $S(U; B_0 ,..., B_n)$ is the set of rational points on $U$ which may be represented by an integral $(n+1)$-tuple $(x_0 ,..., x_n)$ with $|x_m| \leq B_m$ for $m \in \{0 ,…, n\}$. If $B_0 = ... = B_n = B$, then we denote this set by $S(U; B)$.

(ii) $S_1(U; B)$ is the set of rational points on $U$ which may be represented by an integral $(n+1)$-tuple $(1, x_1 ,..., x_n)$ with $|x_m| \leq B$ for $m \in \{1,…, n\}$.

(iii) $S(U; B_0 ,..., B_n ; P_1 ,..., P_t)$ is the subset of points in $S(U; B_0 ,..., B_n)$ which specialise to $P_i$ on $X_{p_i}$ for each $i \in \{1, ..., t\}$. If $B_0 = ... = B_n = B$, then we write $S(U; B; P_1 ,..., P_t)$ for this set.

(iv) If $x_0(P_i) \neq 0$ for $i \in \{1, ..., t\}$, then $S_1(U; B; P_1 ,..., P_t)$ is the subset of points in $S_1(U; B)$ which specialise to $P_i$ on $X_{p_i}$ for each $i \in \{1, ..., t\}$.

(v) $N(U; B_0 ,..., B_n)$ (resp. $N(U; B_0 ,..., B_n ; P_1 ,..., P_t)$) is the cardinality of $S(U; B_0 ,..., B_n)$ (resp. $S(U; B_0 ,..., B_n ; P_1, ..., P_t)$). If $B_0 = ... = B_n = B$, then we denote these numbers by $N(U; B)$ (resp. $N(U; B; P_1 ,..., P_t)$).

(vi) $N_1(U; B)$ is the cardinality of $S_1(U; B)$ and $N_1(U; B ; P_1 ,..., P_t)$ is the cardinality of $S_1(U ; P_1, ..., P_t)$.

**Theorem 3.2.** *Let $X \subset \mathbf{P}^n_{\mathbf{Q}} = \text{Proj}(\mathbf{Q}[x_0 ,..., x_n])$ be a reduced equidimensional closed subscheme of dimension r and degree d and $<$ be a graded monomial ordering on the monomials in $(x_0 ,..., x_n)$. Let $(B_0 ,..., B_n) \in \mathbf{R}^{n+1}_{\geq 1}$, $p_1,..., p_t$ be a sequence of different primes and $\mu_1 ,..., \mu_t$ be a sequence of positive integers satisfying*

(3.3) $$p_1^{(d/\mu_1)^{1/r}} \ldots p_t^{(d/\mu_t)^{1/r}} \geq (B_0^{a_0} \ldots B_n^{a_n})^{(r+1)/r}(B_0 \ldots B_n)^\varepsilon.$$

*Finally, let $P_i$ be an $\mathbf{F}_{p_i}$-point of multiplicity $\mu_i$ on $X_{p_i}$ for each $i \in \{1, \ldots, t\}$.*

*Then there exists a form G of degree $O_{d,n,\varepsilon}(1)$ which vanishes at $S(X; B_0, \ldots, B_n; P_1, \ldots, P_t)$ but not at the generic point of X.*

*Proof.* Let $\xi_1, \ldots, \xi_s$ be primitive $(n+1)$-tuples of integers representing the rational points in $S(X; B_0, \ldots, B_n; P_1, \ldots, P_t)$. Let $k \geq 1$ be the unique integer such that $h(k-1) < s \leq h(k)$. Then by lemma 1.11, we may find forms $F_1, \ldots, F_s$ of the same degree $k = (r!/d)^{1/r} s^{1/r} + O_{d,n}(1)$, such that no non-trivial linear combination of them vanishes at the generic point of X and such that

$$|\det(F_j(\xi_l))| \leq s!(B_0^{a_0} \ldots B_n^{a_n})^{kh(k)}(B_0 \ldots B_n)^\kappa,$$

for some $\kappa = O_{d,n}(h(k))$.

Also, by the proof of 1.11 we have that $h(k) = dk^r/r! + O_{d,n}(k^{r-1})$ and $s = dk^r/r! + O_{d,n}(k^{r-1})$. Hence,

$$kh(k) = ks + O_{d,n}(k^r) = (r!/d)^{1/r} s^{1+1/r} + O_{d,n}(s).$$

Thus, if we write $V = B_0 \ldots B_n$ and $W = (B_0^{a_0} \ldots B_n^{a_n})^{(r+1)/r}$, then we get

(3.4) $\log |\det(F_j(\xi_l))| \leq ((r!/d)^{1/r}(r/(r+1)))s \, [s^{1/r}\log W + O_{d,n}(\log V) + O_{d,n}(\log s)]$.

By the *p*-adic estimates in 2.6, there exists a non-negative integer

$$N_i = (r!/\mu_i)^{1/r}(r/(r+1))s^{1+1/r} + O_{d,n}(s),$$

where $p_i^{N_i} | \det(F_j(\xi_l))$.

There exists thus a positive divisor M of $|\det(F_j(\xi_l))|$ such that

$$\log M = (r!/d)^{1/r}(r/(r+1))s^{1+1/r}[(d/\mu_1)^{1/r}\log p_1 + \ldots + (d/\mu_t)^{1/r}\log p_t] + O_{d,n}(s\Sigma_i \log p_i).$$

Set $A = p_1^{(d/\mu_1)^{1/r}} \ldots p_t^{(d/\mu_t)^{1/r}}$. Then,

(3.5) $\log M = (r!/d)^{1/r}(r/(r+1))s(s^{1/r}\log A + O_{d,n}(\log A))$.

Further, by (3.4) and (3.5) we get that $\log M > \log |\det(F_j(\xi_l))|$ as soon as

(3.6) $s^{1/r}\log A + O_{d,n}(\log A) > s^{1/r}\log W + O_{d,n}(\log V) + O_{d,n}(\log s).$

Suppose now that ε>0 is given. Then there exists a positive integer $s_0$ depending solely on $d$, $n$ and ε such that (3.6) holds for all $s \geq s_0$ and all $A \geq V^\epsilon W$. Hence, if $s \geq s_0$ and $A \geq V^\epsilon W$, then $|\det(F_j(\xi_l))| = 0$.

Let $s = s_0 = O_{d,n,\epsilon}(1)$, $k = (r!/d)^{1/r} s^{1/r} + O_{d,n}(1)$ be as above and $\xi_1, ..., \xi_S$ be a set of primitive $(n+1)$-tuples of integers representing all rational points in $S(X; B_0, ..., B_n; P_1, ..., P_t)$. Then we conclude from the above that there are forms $F_1, ..., F_s$ of degree $k$ such that no non-trivial linear combination of them vanishes at the generic point of $X$ and such that the $s \times S$-matrix $F_j(\xi_l)$ where $1 \leq j \leq s$ and $1 \leq l \leq S$ is of rank less than $s$. There exists therefore a linear combination $G$ of $F_1, ..., F_s$ with the desired properties. This completes the proof.

**Corollary 3.7**. *Let $X \subset \mathbf{P}^n_\mathbf{Q} = \mathrm{Proj}(\mathbf{Q}[x_0, ..., x_n])$ be a reduced equidimensional closed subscheme of dimension r and degree d such that the hyperplane $\Pi_0 \subset \mathbf{P}^n$ defined by $x_0 = 0$ intersects X properly. Let $B \geq 1$, $p_1, ..., p_t$ be a sequence of different primes and $\mu_1, ..., \mu_t$ be a sequence of positive integers satisfying*

(3.8) $\qquad\qquad p_1^{(d/\mu_1)^{1/r}} \cdots p_t^{(d/\mu_t)^{1/r}} \geq B^{1+\epsilon}.$

*Finally, let $P_i$ be an $\mathbf{F}_{p_i}$-point of multiplicity $\mu_i$ on $X_{p_i}$ with $x_0(P_i) \neq 0$ for each $i \in \{1, ..., t\}$.*

*Then there exists a form G of degree $O_{d,n,\epsilon}(1)$ which vanishes at $S_1(X; B; P_1, ..., P_t)$ but not at the generic point of X.*

*Proof.* This follows theorem 3.2 applied to the case where $(B_0, ..., B_n) = (1, B, ..., B)$ and to a graded monomial ordering as in lemma 1.12.

**Remark 3.9**. In the applications $X$ will be the closure of a closed subscheme $X_0 \subset \mathbf{A}^n_\mathbf{Q}$ under the open embedding of $\mathbf{A}^n_\mathbf{Q} = \mathrm{Spec}(\mathbf{Q}[x_1/x_0, ..., x_n/x_0])$ in $\mathbf{P}^n_\mathbf{Q} = \mathrm{Proj}(\mathbf{Q}[x_0, ..., x_n])$. Then $\Pi_0 \subset \mathbf{P}^n$ intersects $X$ properly.

## 4. Rational points on space curves

We shall in this section study $N_1(C; B)$ for space curves. If $F(x_0, ..., x_n)$ is a form of positive degree $d$ with rational coefficients, then we may regard it as a linear form in the monomials of degree $d$ in $(x_0, ..., x_n)$. We will denote by $H(F)$ the projective height of this linear form.

**Lemma 4.1**. *Let $C \subset \mathbf{P}^3$ be an integral curve of degree d over $\mathbf{Q}$ which is not geometrically integral. Then $\#C(\mathbf{Q}) = O_d(1)$.*

*Proof.* This follows from the arguments in the proof of [Sa$_1$], lemma 2.1.

**Lemma 4.2**. *Let $C \subset \mathbf{P}^3$ be a closed subscheme defined over $\mathbf{Q}$ where all irreducible components are of dimension at most one. Then*

$$N_1(C; B) = O_D(B),$$

where $D$ is the sum of the degrees of the irreducible components of $C$.

*Proof.* Let $C_a \subset C$ be the closed subscheme defined by the equation $x_3 = ax_0$ and $C_{a,b} \subset C_a$ be the closed subscheme defined by the equation $x_2 = bx_0$. Then,

$$N_1(C\,;B) = \sum_{a=-B}^{B} N_1(C_a;B), \qquad N_1(C_a\,;B) = \sum_{b=-B}^{B} N_1(C_{a,b};B).$$

In the first sum, we have $N_1(C_a\,;B) \leq D$ as soon as $C_a$ is of dimension 0. This happens for all but $O_D(1)$ values of $a$. For the remaining $a$, we use that $N_1(C_{a,b}\,;B) \leq D$ for all but $O_D(1)$ values of $b$.

**Lemma 4.3**. *Let $C \subset \mathbf{P}^3$ be an integral curve of degree $d$ defined over $\mathbf{Q}$ and $C^{cl}$ be the scheme-theoretic closure of $C$ in $\mathbf{P}^3_{\mathbf{Z}}$. Let $P$ be an $\mathbf{F}_p$-point on $C_p = C^{cl} \times_{\mathbf{Z}} \mathbf{F}_p$ of multiplicity $\mu$ with $x_0(P) \neq 0$. Then*

$$N_1(C\,;B;P) = O_{d,\varepsilon}(B^{1/d+\varepsilon}/p^{1/\mu} + B^\varepsilon).$$

*Proof.* It suffices by lemma 4.1 to consider the case where $C$ is geometrically integral. From lemma 1.3 and lemma 1.4 we conclude that the scheme $C$ is defined by forms $F_1, \ldots, F_t$ of degree $\ll_d 1$ and we may assume that these forms are linearly independent so that $t \ll_d 1$. By [Br], lemma 5 then either there exists a form $G$ of degree $\ll_d 1$ which vanishes at $S(C\,;B)$ but not at the generic point of $C$ or there exists $k \ll_d 1$ such that :

(4.4) $\qquad H(F_i) = O_d(B^k)$, for $i \in \{1, \ldots, t\}$.

In the first case we obtain that $N(C;B) \leq d(\deg G) \ll_d 1$ by the theorem of Bézout (see [Fu$_1$], Ex. 8.4.6.) We may and shall therefore assume that (4.4) holds.

By corollary 3.7 there exists for each prime $q \geq B^{(1/d)(1+\varepsilon)}/p^{1/\mu}$ and each non-singular $\mathbf{F}_q$-point $Q$ on $C_q$ a form $G$ of degree $O_{d,\varepsilon}(1)$ which vanishes at $S_1(C;B;P,Q)$ but not at the generic point of $C$. Hence $N_1(C;B;P,Q) \leq d(\deg G) \ll_{d,\varepsilon} 1$ by the theorem of Bézout.

If $\boldsymbol{\beta} = (\beta_0, \beta_1, \beta_2, \beta_3) \in \mathbf{Z}^4$ represents a non-singular point on $C$, then there are two forms $F, G \in \{F_1, \ldots, F_t\}$ and $0 \leq i < j \leq 3$ such that $\Phi = (\partial F/\partial x_i)(\partial G/\partial x_j) - (\partial F/\partial x_j)(\partial G/\partial x_i)$ does not vanish at $\boldsymbol{\beta}$. Moreover, if $\boldsymbol{\beta}$ represents a point in $S_1(C;B;P)$, then there exists $l \ll_d 1$ with $|\Phi(\boldsymbol{\beta})| \ll O_d(B^l)$. There is therefore a positive integer $N_0 = O_{d,\varepsilon}(1)$ such that the product of any set of $N \geq N_0$ primes $q \geq B^\varepsilon$ primes will exceed $|\Phi(\boldsymbol{\beta})|$. By Bertrand's postulate we may thus find a set $\Omega$ of $O_{d,\varepsilon}(1)$ primes $q$ with the following properties.

(4.5) $B^{(1/d)(1+\varepsilon)}/p^{1/\mu} + B^\varepsilon \leq q \ll_{d,\varepsilon} B^{(1/d)(1+\varepsilon)}/p^{1/\mu} + B^\varepsilon \qquad$ for all $q \in \Omega$

(4.6) Any non-singular point in $S_1(C;B;P)$ specialises to a non-singular $\mathbf{F}_q$-point on $C_q$ for some $q \in \Omega$.

By summing over the set of all non-singular $\mathbf{F}_q$-points $Q$ on $C_q$ with $x_0(Q) \neq 0$ for all $q \in \Omega$, we obtain :

$N_1(C; B; P) \leq \Sigma_Q N_1(C; B; P,Q) <<_{d,\varepsilon} \Sigma_Q 1 <<_{d,\varepsilon} (\#\Omega)(B^{(1/d)(1+\varepsilon)}/p^{1/\mu} + B^{\varepsilon}) << B^{(1/d)(1+\varepsilon)}/p^{1/\mu} + B^{\varepsilon}$.

This completes the proof.

## 5. Universal polynomials and the geometry of surfaces in $\mathbf{P}^3$

We shall in this section study certain geometric conditions (%$_1$)-(%$_4$) for points on surfaces in $\mathbf{P}^3$. These conditions define an open subscheme of the surface which is non-empty for any surface of degree $d \geq 4$ over a field of characteristic 0. We also give an upper bound (see 5.6) for the multiplicity of divisors on the surface at points on this open set. This bound will be important in the next section when we apply the determinant method.

**Notation 5.1**. Let $\underline{k}=(k_0,..., k_n)$ be an $(n+1)$-tuple of non-negative integers. Then $deg(\underline{k}) = k_0 + ... + k_n$. Also, if $x=(x_0,..., x_n)$, then $x^{\underline{k}} = x_0^{k_0} ... x_n^{k_n}$.

**Definition 5.2**. Let $X \subset \mathbf{P}^3$ be the surface defined by a form $F(x_0,..., x_3)$ with coefficients in an algebraically closed field K. Then we denote by (%$_j$), $j \in \{1,2,3,4,5\}$, the following geometric conditions for a closed point P on X.

(%$_1$) P is a non-singular point.
(%$_2$) P does not belong to any plane contained in X and P is of multiplicity at most two on any plane section of X.
(%$_3$) P does not lie on any line on X.
(%$_4$) P is a non-singular point on any plane cubic contained in X.
(%$_5$) P satisfies all the previous conditions (%$_1$) - (%$_4$).

We do not demand that the plane cubics which occur in this section are irreducible or reduced. Note also that it suffices to consider the tangent plane $T_{X,P}$ in (%$_2$) and (%$_4$) if (%$_1$) holds since any other plane intersect X transversally at P. This gives another equivalent way of stating (%$_5$), which we shall use in the proof of theorem 6.4.

**Lemma 5.3**. Let d be a positive integer and let $\underline{k}=(k_0,..., k_3)$ run over all quadruples of non-negative integers with $deg(\underline{m}) = d$. Then there exists for each d and each $j \in \{1,2,3,4,5\}$ a finite set of universal bihomogeneous polynomials $H_0(a_{\underline{k}}; x_0,..., x_3), ..., H_t(a_{\underline{k}}; x_0,..., x_3)$ with integer coefficients such that the following assertion holds for any algebraically closed field K and any form

$$F(x_0,..., x_3) = \sum a_{\underline{k}} x^{\underline{k}}, \qquad a_{\underline{k}} \in K$$

of degree d.

Let $X = \mathrm{Proj}\, K[x_0,..., x_3]/(F(x_0,..., x_3))$ and P be a closed point on X represented by a quadruple $(\beta_0,..., \beta_3) \neq (0,..., 0) \in K^4$. Then satisfies (%$_j$) if and only if $H_i(a_{\underline{k}}; \beta_0,..., \beta_3) \neq 0$ for some $i \in \{0,..., t\}$.

*Proof.* (%$_1$) We apply the Jacobian criterion and choose $H_i = G_{x_i}$ for $i \in \{0,..., n\}$.

($\%_2$) Let $\mathbf{P}^{3\vee}$ be the Hilbert scheme of all planes $\Pi \subset \mathbf{P}^3$ and $H$ the Hilbert scheme of all surfaces $X \subset \mathbf{P}^3$ of degree $d$ in $\mathbf{P}^3$. Then there exists a closed subscheme $I_1 \subset \mathbf{P}^{3\vee} \times H \times \mathbf{P}^3$ representing triples $(\Pi, X, P)$ where $\Pi \cap X$ has intersection multiplicity $\geq 3$ at $P$. (If $P \in \Pi$ and $\Pi \subseteq X$, then we formally define the intersection multiplicity of $\Pi \cap X$ at $P$ to be $+\infty$.) Now let $Z_1$ be the scheme-theoretic image of $I_1$ under the projection from $\mathbf{P}^{3\vee} \times H \times \mathbf{P}^3$ to $H \times \mathbf{P}^3$. Then $Z_1$ is a closed subscheme of $H \times \mathbf{P}^3$ by the main theorem in elimination theory. It is defined by finitely many bihomogeneous polynomials in $(a_{\underline{k}}; x_0,..., x_3)$. (For a more explicit construction of these polynomials, see the proof of lemma 11 in [BwHeSa].)

($\%_3$) Let $G$ be the Grasmannian of all lines on $\mathbf{P}^3$ and $H$ the Hilbert scheme of all surfaces $X$ of degree $d$ in $\mathbf{P}^3$. Then the exists a closed subscheme $I_2 \subset G \times H \times \mathbf{P}^3$ representing triples $(L, X, P)$ where $P$ is a point on the line $L$ lying on the surface $X \subset \mathbf{P}^3$ of degree $d$. Now let $Z_2$ be the scheme-theoretic image of $I_2$ under the projection from $G \times H \times \mathbf{P}^3$ to $H \times \mathbf{P}^3$. It is defined by finitely many bihomogeneous polynomials in $(a_{\underline{k}}; x_0, ... , x_3)$.

($\%_4$) Let $H'$ be the Hilbert scheme of all cubics $C \subset \mathbf{P}^3$ contained in a plane and $H$ the Hilbert scheme of all surfaces $X$ of degree $d$ in $\mathbf{P}^3$. Then there is a closed subscheme $I_3 \subset H' \times H \times \mathbf{P}^3$ representing triples $(C, X, P)$ where $C \subset X$ and $C$ has intersection multiplicity $\geq 2$ at $P$. Now let $Z_3$ be the scheme-theoretic image of $I_3$ under the projection from $G \times H \times \mathbf{P}^3$ to $H \times \mathbf{P}^3$. It is defined by finitely many bihomogeneous polynomials in $(a_{\underline{k}}; x_0,..., x_3)$.

($\%_5$) It suffices to take products of the universal polynomials that occur in ($\%_1$)-($\%_4$).

To prove the next lemma, we shall need the theorem of Kronecker-Castelnuovo [Ca]. It says (cf. p.130 in [SR]) that a surface with a two-dimensional family of reducible hyperplane sections is either ruled by lines or is a linear projection of the Veronese surface. For a modern rigorous version of the proof of Castelnuovo, see [Rg$_1$].

**Lemma 5.4**. *Let $X \subset \mathbf{P}^3$ be an integral surface of degree $d \geq 4$ over an algebraically closed field of characteristic 0. Suppose that there are only finitely many lines on $X$. Then the closed subset $W \subseteq \mathbf{P}^{3\vee}$ of planes $\Pi$ such that $\Pi \cap X$ contains a cubic curve is of dimension $\leq 1$.*

*Proof.* If dim $W \geq 2$, then there exists by the theorem of Kronecker and Castelnuovo a linear projection of the Veronese surface $V_4 \subset \mathbf{P}^5$ onto $X$. Hence $d = \deg X \leq \deg V_4 = 4$ (cf. p.76 in [Mu$_1$] and [IP], 4.5.5) so that $d = 4$ by assumption. There is therefore a residual line on any any plane section of $X$ which contains a cubic. But then dim $W \leq 1$ as there are only finitely many lines on $X$.

The quartic surfaces in $\mathbf{P}^3$, which appear in the proof are called Steiner surfaces. It is known that they contain three double lines (cf. [IP], 4.5.5 and pp.133-134 in [SR]).

**Lemma 5.5**. *Let $X \subset \mathbf{P}^3$ be an integral surface of degree $d$ over an algebraically closed field of characteristic 0. Suppose that there are only finitely many lines on X. Then there exists a closed point $P$ on $X$ satisfying ($\%_1$)-($\%_3$) if $d \geq 3$ and satisfying ($\%_1$)-($\%_4$) if $d \geq 4$.*

*Proof.* It follows from 5.3 that the set $U_i$ of points satisfying ($\%_i$) is open in the Zariski topology. It is therefore enough to show that there is a closed point $P_i$ on $U_i$ satisfying ($\%_i$) for $i=1,2,3,4$. This is obvious for ($\%_1$) and ($\%_3$) and it is shown in [BwHeSa], lemma 11 for ($\%_2$). It is thus enough to show that there exists a closed point $P=P_4$ satisfying ($\%_4$) when $d\geq 4$.

Let $\gamma: X_{ns} \to \mathbf{P}^{3\vee}$ be the Gauss map which sends a non-singular closed point $P$ to the tangent plane $T_{X,P}$ of $X$ at $P$. Then, $\dim \gamma(X_{ns})=2$ by [Z], thm I.2.3c as $X$ is not ruled by lines. Hence by lemma 5.4, we may find a non-singular closed point $P$ of $X$ such that $T_{X,P} \cap X$ contains no cubic curve. This implies that $P$ satisfies ($\%_4$) since any plane curve containing $P$ as a singular point must be contained in its tangent plane.

**Lemma 5.6.** *Let $X \subset \mathbf{P}^3$ be an integral surface of degree$\geq 4$ over an algebraically closed field and $P$ be a closed point satisfying ($\%_1$) - ($\%_4$). Let $C$ be a closed purely one-dimensional closed subscheme of $X$ of degree $e$ and $\mu$ be the multiplicity of $C$ at $P$. Then $\mu \leq e/2$.*

*Proof.* We may assume that $C$ is integral because of the remark in [Fu$_1$], Ex. 4.3.4. By ($\%_1$), $P$ is a non-singular point on $X$. Let $T_{X,P}$ be tangent plane of $X$ at $P$. Then, if $C \subset T_{X,P}$, we deduce from ($\%_2$) that $\mu \leq 2$ and from ($\%_4$) that $\mu \leq 1$ in case $e=3$. Also, if $e=2$, then $\mu \leq 1$ since any geometrically integral conic is non-singular. Finally, if $C$ is a line, then $\mu=0$ by ($\%_3$). We have thus shown that $\mu \leq e/2$ for integral curves on $T_{X,P}$. If $C \not\subset T_{X,P}$, then we apply the following lemma.

**Lemma 5.7.** *Let $X \subset \mathbf{P}^3$ be an integral surface over an algebraically closed field. Let $P$ be a non-singular closed point on $X$ and $C$ be an integral curve on $X$ passing through $P$ which is not contained in the tangent plane $T_{X,P}$ of $X$ at $P$. Then the multiplicity of $C$ at $P$ is at most half the degree of $C$.*

*Proof.* Let $D$ be the scheme-theoretic intersection $T_{X,P} \cap X$. It is a purely one-dimensional subscheme of $X$ which intersects $C$ properly at $P$. Let $i(P, C \cdot D; X)$ be the local intersection multiplicity of $C$ and $D$ at $P$ and $\mu_P(C)$ resp. $\mu_P(D)$ be the multiplicities of $C$ resp. $D$ at $P$. Then, $i(P, C \cdot D; X) \geq \mu_P(C)\mu_P(D)$ by a result of M. Noether (see [Fu$_1$], Ex.7.1.11, [Fu$_2$], 2.8) or by [Fu$_1$], Cor.12.4 applied to an open non-singular neighbourhood of $P$. Also, $i(P, C \cdot D; X) = i(P, C \cdot T_{X,P}; \mathbf{P}^3) \leq \deg C$. Therefore, $\mu_P(C) \leq \deg C / \mu_P(D) \leq \deg C / 2$ as $P$ is a singular point on $D = T_{X,P} \cap X$. This completes the proof.

## 6. Integral points on surfaces in $\mathbf{A}^3$

The aim of this section is to prove theorem 0.4 for polynomials $f(y_1, y_2, y_3) \in \mathbf{Z}[y_1, y_2, y_3]$. We shall thereby compactify the affine scheme $Y \subset \mathbf{A}^3_{\mathbf{Z}}$ defined by $f$ to a projective scheme $\Xi \subset \mathbf{P}^3_{\mathbf{Z}}$. We use the results in section 5 to prove that the rational points on $X = \Xi_{\mathbf{Q}}$ which satisfy ($\%_1$)-($\%_4$) will specialise to $\mathbf{F}_p$-points of bounded multiplicity for suitably chosen primes $p$ if the coefficients of $f$ are not too big. These multiplicity estimates are crucial in order to make efficient use of lemma 4.3 when counting points of bounded height on curves on $Y$.

**Notation 6.1.** Let $Y \subset \mathbf{A}^n_{\mathbf{Z}}$ be the affine scheme defined by a polynomial $f(y_1,\ldots, y_n)$ with integer coefficients. Then $n(Y; B)$ is the number of $n$-tuples $\mathbf{y} = (y_1,\ldots, y_n)$ in integers such that $y_1,\ldots, y_n \in [-B, B]$ and $f(\mathbf{y}) = 0$.

**Lemma 6.2**. *Let $X \subset \mathbf{P}_\mathbf{Q}^3$ be an integral surface over $\mathbf{Q}$ of degree $d$ such that the hyperplane $\Pi_0 \subset \mathbf{P}^3$ defined by $x_0 = 0$ intersects $X$ properly. Let $p$ be a prime such that $p \geq B^{1/\sqrt{d}+\varepsilon}$ and $P$ be a non-singular $\mathbf{F}_p$-point on $X_p$. Then there exists a form $G$ of degree bounded in terms of $d$ and $\varepsilon$ which vanishes at $S_1(X; B; P)$ but not at the generic point of $X$.*

*Proof.* This is a special case of corollary 3.7 essentially due to Heath-Brown (see [He$_3$], th.15).

**Lemma 6.3**. *Let $F(x_0, x_1, x_2, x_3)$ be an irreducible quaternary form of degree $d \geq 2$ with rational coefficients and $X = \mathrm{Proj}(\mathbf{Q}[x_0, x_1, x_2, x_3]/(F))$. Then either there exists a form $G(x_0, x_1, x_2, x_3) \in \mathbf{Q}[x_0, x_1, x_2, x_3]$ of degree $d$ not proportional to $F$ which vanishes on $S(X; B)$ or $H(F) \ll_d B^{d(d+1)(d+2)(d+3)/6}$.*

*Proof.* See [BwHe$_1$], lemma 3 and the proof of [He$_2$], th. 4.

**Theorem 6.4**. *Let $f(y_1, y_2, y_3) \in \mathbf{Z}[y_1, y_2, y_3]$ be a polynomial of degree $d \geq 4$ which is absolutely irreducible over $\mathbf{Q}$. Suppose that there are only finitely many lines over $\overline{\mathbf{Q}}$ on the surface $Y_\mathbf{Q} \subset \mathbf{A}_\mathbf{Q}^3$ defined by $f(\mathbf{y}) = 0$. Then,*

$$n(Y; B) = O_{d,\varepsilon}(B^{1+\varepsilon}).$$

*Proof.* We may and shall assume that the coefficients of $f$ have no common factor. Then $F(x_0, x_1, x_2, x_3) = x_0^d f(x_1/x_0, x_2/x_0, x_3/x_0)$ is a primitive form with integer coefficients and $\Xi = \mathrm{Proj}(\mathbf{Z}[x_0, x_1, x_2, x_3]/(F))$ the scheme-theoretic closure of $X = \mathrm{Proj}(\mathbf{Q}[x_0, x_1, x_2, x_3]/(F))$ in $\mathbf{P}_\mathbf{Z}^3$. Also, $n(Y; B) = N_1(X; B)$ (cf. 3.1) so that $n(Y; B) \leq N(X; B)$.

Let $B > 1$. By lemma 6.3 there are two possibilities. In the first case there exists a form $G(x_0, x_1, x_2, x_3)$ of degree $d$ not proportional to $F$ which vanishes on $S(X; B)$. Let $W$ be the closed subscheme of $X$ defined by $G = 0$ and $D$ be the sum of the degrees of all irreducible components of $W$. Then $n(Y; B) = N_1(X; B) = N_1(W; B) \leq O_D(B)$ by lemma 4.2 and $D \leq d^2$ by the version of Bézout's theorem in [Fu$_1$], Ex. 8.4.6. Hence we obtain that $N(X; B) = O_d(B)$, which is satisfactory. We may and shall thus from now on assume that we are in the second case where $H(F) \ll_d B^{d(d+1)(d+2)(d+3)/6}$.

Let $d \geq 4$ and $\underline{k}$ run over all quadruples $(k_0, k_1, k_2, k_3)$ of non-negative integers such that $k_0 + k_1 + k_2 + k_3 = d$. There exists, then, by lemma 5.3 a finite set $T$ of bihomogeneous polynomials $H_i(a_{\underline{k}}; x_0, x_1, x_2, x_3)$, $i = 1, \ldots, t$ such that the following holds for any algebraically closed field $K$ and any form $F(x_0, x_1, x_2, x_3) = \Sigma\, a_{\underline{k}} x^{\underline{k}}$ of degree $d$ in $K[x_0, x_1, x_2, x_3]$.

(6.5) Let $(\beta_0, \beta_1, \beta_2, \beta_3) \neq (0, 0, 0, 0)$ be a quadruple in $K$ representing a closed point $P$ on $X = \mathrm{Proj}(K[x_0, x_1, x_2, x_3]/(F))$. Then $H_i(a_{\underline{k}}; \beta_0, \beta_1, \beta_2, \beta_3) \neq 0$ for some $H_i$ in $T$ if and only if the following conditions hold.

(a) $P$ is a non-singular point on $X$.
(b) The intersection multiplicity of $X$ with the tangent plane $T_{X,P}$ of $X$ at $P$ is at most two.
(c) $P$ does not lie on any line on $X$.
(d) $P$ is a non-singular point on any cubic curve on $X \cap T_{X,P}$ passing through $P$.

Also, let $G_i(x_0, x_1, x_2, x_3) = H_i(\underline{a_k}; x_0, x_1, x_2, x_3)$, $i \in \{1, \ldots, t\}$ be the forms obtained by specialising at the coefficients $\underline{a_k}$ of $F(x_0, x_1, x_2, x_3)$. Let $Z$ be the closed subscheme of $X$ defined by the homogeneous ideal $(F, G_1, \ldots, G_t)$ of $\mathbf{Q}[x_0, x_1, x_2, x_3]$. Then dim $Z \leq 1$ by lemma 5.5. Further, as $t$ and the degrees of the forms $(F, G_1, \ldots, G_t)$ are bounded in terms of $d$, we conclude from [Fu$_1$], Ex. 8.4.6, that the sum of the degrees of the irreducible components is bounded in terms of $d$. Hence $N_1(Z; B) = O_d(B)$ by lemma 4.2.

It remains to prove that $N_1(U; B) = O_{d,\varepsilon}(B^{1+\varepsilon})$ for $U = X \setminus Z$. On applying the same arguments as in the proof of lemma 4.3 we find a set of primes $\Omega$ of $O_{d,\varepsilon}(1)$ primes $p$ such that the following holds.

(6.6)(a) $B^{1/\sqrt{d}+\varepsilon} \leq p <<_{d,\varepsilon} B^{1/\sqrt{d}+\varepsilon}$ for each prime $p \in \Omega$.
(b) there exists for each point $(1, b_1, b_2, b_3)$ in $S_1(U; B)$ a prime $p \in \Omega$ such that $p$ does not divide some $G_i(1, b_1, b_2, b_3)$, $i \in \{1, \ldots, t\}$.

We have thus by (6.5) and (6.6)(b) that $N_1(U; B) \leq \sum_P N_1(X; B; P)$ where $P$ runs over all $\mathbf{F}_p$-points with $x_0(P) \neq 0$ such that $p \in \Omega$ and

(6.7)(a) $P$ is a non-singular $\mathbf{F}_p$-point on $X_p$.
(b) The intersection multiplicity of $X_p$ with the tangent plane $T_{X_p, P}$ of $X_p$ at $P$ is at most two.
(c) $P$ does not lie on any line on $X_p \times K_p$
(d) $P$ is a non-singular point on any cubic curve on $(T_{X_p, P} \cap X_p) \times K_p$ passing through $P$.

Here $K_p$ denotes an algebraic closure of $\mathbf{F}_p$.

As $\# X_p(\mathbf{F}_p) = O_d(p^2)$ and $\#\Omega = O_{d,\varepsilon}(1)$, we see also from (6.6)(a) that there are $<<_{d,\varepsilon} B^{2/\sqrt{d}+\varepsilon}$ such closed points $P$ on $\Xi$. To complete the proof, it is thus enough to show that

(6.8) $\quad N_1(X; B; P) <<_{d,\varepsilon} B^{1-2/\sqrt{d}+\varepsilon} + B^\varepsilon$,

for each such closed point $P$.

To prove (6.8), we note that the hypothesis of lemma 6.2 is satisfied by (6.6)(a). There is thus a form $G$ of degree bounded in terms of $d, \varepsilon$ which vanishes at $S(X; 1, B, B, B; P)$ but not at the generic point of $X$. There is thus by the theorem of Bézout a set of $O_{d,\varepsilon}(1)$ curves $C$ of degree $O_{d,\varepsilon}(1)$ such that $S_1(X; B; P) = \bigcup_C S_1(C; B; P)$. This implies that to prove (6.8), it will suffice to establish the estimate

(6.9) $\quad N_1(C; B; P) <<_{d,\varepsilon} B^{1/e-2/e\sqrt{d}+\varepsilon} + B^\varepsilon$

for curves $C$ of degree $e$ on $X$. But $C_p = C^{cl} \times_\mathbf{Z} \mathbf{F}_p \subset \mathbf{P}^3_{\mathbf{F}_p}$ is of the same degree as $C \subset \mathbf{P}^3_\mathbf{Q}$ since $C^{cl}$ is flat over Spec $\mathbf{Z}$ (see III.9.8 and III.9.9 in [Ha]). The multiplicity of $C_p$ at $P$ is thus at most $e/2$ by (6.7) and lemma 5.6. We therefore obtain (6.9) from lemma 4.3. This completes the proof of theorem 6.4.

# 7. Projective varieties with many linear subspaces

We shall in the next section generalise theorem 6.4 to affine hypersurfaces of dimension $\geq 2$. The proof will be inductive and based on an analysis of the hyperplane sections of a projective closure of the affine hypersurface. For the induction to work we must relate the geometry of the variety to the geometry of its hyperplane sections. We shall therefore in this section give a criterion 7.5 for when the union of all divisors of degree one contains all closed points on a projective hypersurface.

**Lemma 7.1**. *Let k be an algebraically closed field of characteristic $0$ and $Z \subset \mathbf{P}^r$ be a (possibly non-reduced or reducible) hypersurface over k of dimension at least two. Let P be a closed point on Z such that there is an (r-2)-plane on $\Pi \cap Z$ through P, for any hyperplane $\Pi \subset \mathbf{P}^r$ through P. If r=3, suppose also that there is a plane $\Pi$ through P, such that there is only one line on $\Pi \cap Z$ through P and that this line is simple on $\Pi \cap Z$. Then there is an (r-1)-plane on Z passing through P.*

*Proof.* Let $\Gamma$ be the hyperplane in the dual projective space $\mathbf{P}^{r\vee}$ parameterising hyperplanes $\Pi \subset \mathbf{P}^r$ through $P$ and $F_{r-2}(Z; P)$ be the Hilbert scheme of $(r-2)$-planes $\Lambda$ on $Z$ through $P$. Let $I \subset F_{r-2}(Z; P) \times \Gamma$ be the closed subscheme representing pairs $(\Lambda, \Pi)$ where $\Lambda \subset \Pi$ and let $p: I \to \Gamma$ be the second projection. Then $\Gamma = p(I)$ by assumption. Next, let $Z_1, \ldots, Z_m$ be the irreducible components of $Z$ containing $P$ with their reduced scheme structures and let $I_j \subset F_{r-2}(Z_j; P) \times \Gamma$ and $p_j: I_j \to \Gamma$, $1 \leq j \leq m$ be defined in the same way as $I$ and $p$. Then $p_j(I_j)$ is a closed subset of $\Gamma$ for $1 \leq j \leq m$ since $F_{r-2}(Z_j; P)$ is a proper $k$-scheme. Also, $p(I) = \cup_j p_j(I_j)$ since any $(r-2)$-plane $\Lambda$ on $Z$ which passes through $P$ is contained in one of the components $Z_j$, $j \in \{1, \ldots, m\}$. Hence, as $\Gamma = \cup_j p_j(I_j)$ is irreducible, we must have that $p_j(I_j) = \Gamma$ for one of the irreducible components $Z_j$ of $Z$. There exists thus an irreducible component $Z_j$ of $Z$ containing $P$ such that any hyperplane section $\Pi \cap Z_j$ through $P$ contains an $(r-2)$-plane through $P$. Let $q_j: Z_j - P \to \Pi' = \mathbf{P}^{r-1}$ be the projection from $P$ to a hyperplane $\Pi' \subset \mathbf{P}^r$ not containing $P$. If $\dim q_j(Z_j - P) = 1$, then $r = 3$ and $Z_j$ is a two-dimensional cone with vertex $P$. Hence $Z_j$ is a plane since otherwise we would have two different lines or a line of multiplicity $\geq 2$ through $P$ for any plane section of $Z_j$. If $\dim q_j(Z_j - P) \geq 2$, then we apply Bertini´s theorem (cf. [FuL], th. 2.1 and [J], cor. 6.11.3) and conclude that $\Pi \cap Z_j$ is integral and of codimension two in $\mathbf{P}^r$ for a generic hyperplane $\Pi \subset \mathbf{P}^r$ through $P$. But then $\Pi \cap Z_j$ is an $(r-2)$-plane since it contains an $(r-2)$-plane passing through $P$. Hence $Z_j$ is an $(r-1)$-plane on $Z$ passing through $P$, thereby completing the proof.

**Lemma 7.2**. *Let k be an algebraically closed field of characteristic $0$ and $Z \subset \mathbf{P}^3$ be a non-linear integral surface such that any point on Z lies on at least two different lines or on a line in the singular locus of Z. Then $Z \subset \mathbf{P}^3$ is a non-singular quadric.*

*Proof.* Let $U \subseteq Z$ be the set of non-singular points $P$ on $Z$ such that $P$ is of multiplicity two on the intersection of $Z$ with the tangent plane $T_P$ at $P$. Then $U$ is an open non-empty subset of $Z$ (cf. e.g. lemma 10 and 11 in [BwHeSa]). It follows from the assumption that there are exactly two different two lines through points on $U$. There is, therefore, a morphism $\eta : U \to H$ to the Hilbert scheme $H$ of conics in $\mathbf{P}^3$ which sends $P \in U$ to the conic consisting of the two lines on $X$ through $P$. The image $\eta(U)$ lies in the closed subscheme $H_Z \subset H$ of conics on $Z$. Let $\pi : H \to \mathbf{P}^{3\vee}$ be the morphism which maps a conic to the plane it spans (cf. [H$_1$, 1.b]). Then $\gamma = \pi\eta$ is the Gauss map which sends $P \in U$ to its tangent plane $T_P \in \mathbf{P}^{3\vee}$. Also, since $P$ is an ordinary double point on $T_P \cap Z$ for $P \in U$, we conclude from (17) and (15) in chapter I of [Kl$_2$] that the restriction of $\gamma$ to some open subset $V \neq \emptyset$ is injective.

Now fix a (closed) point $P_1$ on $V$ and let $L_1$ be one of the lines on $X$ through $P_1$. Then since $H$ is proper (see [Ha], th,II.4.7) there exists a (unique) morphism $\eta_1: L_1 \to H$ with $\eta_1(L_1) \subseteq H_Z$ that extends the restriction of $\eta$ to $U \cap L_1$. Let $L_1^{\vee} \subset \mathbf{P}^{3\vee}$ be the dual line parameterising planes in $\mathbf{P}^3$ through $L$ and $\gamma_1 = \pi \eta_1$. Then $\gamma_1(U \cap L_1) \subseteq L_1^{\vee}$ since the tangent planes of points on $U \cap L_1$ contain $L_1$. Now as $\gamma_1$ is injective on $V \cap L_1$, we get a birational morphism from $L_1$ to $L_1^{\vee}$ which must be an isomorphism since $L_i \cong L_i^{\vee} \cong \mathbf{P}^1$. Therefore, $\gamma_1$ is a closed immersion. The conic parameterised by $\eta_1(P)$ for a point $P$ on $L_1$ is a singular conic $C \subset Z$ which contains $L_1$. Also, $P_1$ can only be a double point on $C$ if $C$ spans the tangent plane of $Z$ at $P_1$. Hence $C = L \cup L_1$ for a line $L \subset Z$ with $L \cap L_1 = P$. There is thus a morphism $\lambda_1 : L_1 \to F_1(Z)$ which sends $P \in L_1$ to the Hilbert point of $L$ and where $\gamma_1(P)$ parameterises the plane $<L, L_1>$. Now as $\gamma_1$ is injective, the lines parameterised by points on $\lambda_1(L_1)$ must be pairwise disjoint. Hence, if $Z_1 \subset Z \times \lambda_1(L_1)$ is the family of lines induced by the universal family of lines over $F_1(Z)$, then the projection $\mathrm{pr}_1: Z_1 \to Z$ is bijective. Therefore, $\mathrm{pr}_1$ must be an isomorphism since it is proper (see [G$_2$], Prop. 8.12.5). In particular, $Z$ is non-singular and $L^2 = 0$ for any line $L$ in the family parameterised by $\lambda_1(L_1)$. But it follows from the adjunction formula that $L^2 = 2-d$ for a line on a non-singular surface of degree $d$. Hence $Z \subset \mathbf{P}^3$ is a non-singular quadric, as was to be proved.

**Definition 7.3**. Let $k$ be a field and $X$ be a closed equidimensional subscheme of $\mathbf{P}_k^n$ of dimension $r \geq 1$. Let $F_{r-1}(X)$ be the Hilbert scheme of $(r-1)$-planes on $X \subset \mathbf{P}_k^n$ and $I \subset X \times F_{r-1}(X)$ be the universal family of all $(r-1)$-planes on $X$. Then $X$ is said to be *covered by its linear divisors* if the projection from $I \subset X \times F_{r-1}(X)$ to $X$ is surjective.

It is easy to see that $X$ is covered by its linear divisors if and only if each closed point on $\overline{X}$ lies on an $(r-1)$-plane over $\overline{k}$ on $\overline{X}$. We shall also use the following equivalent conditions.

**Proposition 7.4**. *Let $k$ be a field and $X \subset \mathbf{P}_k^n$ be a closed integral $r$-dimensional subscheme. Let $\overline{k}$ be an algebraic closure of $k$ and $\overline{X} = X \times_k \overline{k}$. Then the following conditions are equivalent.*

(a) *$X$ is covered by its linear divisors.*
(b) *$F_{r-1}(X)$ is of positive dimension.*
(c) *There are infinitely many $(r-1)$-planes over $\overline{k}$ on $\overline{X}$.*

*Proof.* (a)$\Rightarrow$(b) If (a) holds, then the projection from $I \subset X \times F_{r-1}(X)$ to $X$ is surjective and dim $I \geq r$. The universal family $I \to F_{r-1}(X)$ of $(r-1)$-planes on $X \subset \mathbf{P}_k^n$ is a $\mathbf{P}^{r-1}$-bundle over $F_{r-1}(X)$. Hence dim $F_{r-1}(X) = $ dim $I - (r-1) \geq 1$.
(b) $\Rightarrow$(c) Suppose there were only finitely many $(r-1)$-planes over $\overline{k}$ on $\overline{X}$. Then there would be only finitely many $\overline{k}$-points on $F_{r-1}(X)$ and hence dim $F_{r-1}(X) \leq 0$.
(c) $\Rightarrow$(a) Let $Z$ be the scheme-theoretic image of $I$ (cf. (7.3)) under the projection from $X \times F_{r-1}(X)$ to $X$. Then, if (a) were false, we would have dim $Z \leq r-1$ by the irreducibility of $X$. Also, any $(r-1)$-plane over $\overline{k}$ on $\overline{X}$ must then be one of the finitely many irreducible components of $\overline{Z} = Z \times_k \overline{k}$. This contradiction completes the proof.

**Theorem 7.5**. *Let k be an algebraically closed field of characteristic* 0 *and* $X \subset \mathbf{P}^n$ *be an integral non-linear hypersurface over k of dimension at least four or of dimension three and degree at least three. Suppose that all hyperplane sections of* $X \subset \mathbf{P}^n$ *are covered by their linear divisors. Then X is covered by its linear divisors.*

*Proof.* Let $I \subset X \times F_{n-2}(X)$ be the universal family of $(n-2)$-planes on $X$ and $\text{pr}_1 : I \to X$ be the projection of $I \subset X \times F_{n-2}(X)$ onto $X$. Then $\text{pr}_1(I)$ is a closed subset of $X$. It is therefore enough to prove that any point $P$ on a non-empty open subset of $X$ lies on an $(n-2)$-plane on $X$. Let $T_P \subset \mathbf{P}^n$ be the tangent plane of $X$ at a non-singular point $P$, $Z = T_P \cap X$ and $\Pi$ be a hyperplane in $T_P$ passing through $P$. Then $\Pi = T_P \cap \Pi_0$ for some hyperplane $\Pi_0 \subset \mathbf{P}^n$ and there exists by the hypothesis an $(n-3)$-plane $\Lambda$ on $\Pi_0 \cap X$ passing through $P$. But any $(n-3)$-plane on $X$ which passes through $P$ is contained in $T_P$. There is thus an $(r-2)$-plane on $\Pi \cap Z = T_P \cap \Pi_0 \cap X$, $r = n-1$ through $P$ for any hyperplane $\Pi$ in $T_P = \mathbf{P}^r$. Hence, if $r \geq 4$, we conclude from lemma 7.1 that any non-singular point $P$ on $X$ lies on an $(r-1)$-plane on $Z = T_P \cap X$.

We obtain the same conclusion from lemma 7.1 when $r = 3$ provided that we can show that the additional hypothesis on $Z = T_P \cap X \subset T_P = \mathbf{P}^3$ in lemma 7.1 holds. We shall do this for non-singular points $P$ outside a closed subset $R$ of $X$. This set $R \subset X$ will be the set of points $P \in X$ where $\dim p^{-1}(P) \geq 3$ for the projection $p$ from a suitable closed subset $D \subset \mathbf{P}^{4\vee} \times X$ to $X$.

To define $D$ and $R$, let $\Sigma \subset \mathbf{P}^{4\vee} \times X$ be the incidence variety representing pairs $(\Pi, P)$ where $P \in \Pi \cap X$ and $J \subset \mathbf{P}^{4\vee} \times X \times F_1(X)$ be the closed subvariety representing triples $(\Pi, P, \Lambda)$ where $P \in \Lambda \subset \Pi \cap X$. Let $\pi : J \to \Sigma$ be the morphism obtained by restricting the projection map $\mathbf{P}^{4\vee} \times X \times F_1(X) \to \mathbf{P}^{4\vee} \times X$ and $S \subset \Sigma$ be the open subset of points $s \in \Sigma$ for which $\pi^{-1}(s)$ is of dimension zero. The stalk $(\pi_*O_J)_s$ of the coherent sheaf $\pi_*O_J$ at $s \in S$ is a finitely generated $O_{S,s}$-module and the tensor product $(\pi_*O_J)_s \otimes k(s)$ with the residue field $k(s)$ of $O_{S,s}$ a finitely generated vector space over $k(s)$. Let $\phi(s)$ be the dimension of this vector space over $k(s)$. Then $\phi: S \to \mathbf{N}$ is an upper semicontinuous function by exercise II.5.8(a) in [Ha]. We may even extend it to an upper semicontinuous function $\phi: \Sigma \to \mathbf{N} \cup \{\infty\}$ by letting $\phi(s) = \infty$ when $s \in \Sigma \setminus S$. Then the set $D$ of points $s \in \Sigma$ for which $\phi(s) \geq 2$ is a closed subset of $\Sigma$ and $\mathbf{P}^{4\vee} \times X$. By upper semicontinuity of fibre dimension [G$_2$], Corollaire 13.1.5, we then conclude that $R$ is a closed subset of $X$.

Let $P \in X \setminus R$ and $H_P \subset \mathbf{P}^{4\vee}$ be the hyperplane parameterising hyperplanes $\Pi \subset \mathbf{P}^4$ through $P$. There exists then by the definition of $R$ and the upper semicontinuity of $\phi$ a non-empty open set $U_P \subset H_P$ of hyperplanes $\Pi_0 \subset \mathbf{P}^4$ containing $P$ such that $\phi(\Pi_0, P) \leq 1$. There is thus for such $\Pi_0$ at most one line on $\Pi_0 \cap X$ through $P$ and this line is simple on $\Pi_0 \cap X$. If $P$ is non-singular, choose $\Pi_0 \neq T_P$ and let $\Pi = T_P \cap \Pi_0$. Then $\Pi \subset T_P = \mathbf{P}^3$ satisfies the last hypothesis in lemma 7.1 with respect to $Z = T_P \cap X$. We obtain therefore from lemma 7.1 that any non-singular point $P$ on $X \setminus R$ lies on a plane on $Z = T_P \cap X$.

It remains to prove that $R \neq X$. Suppose instead that $R = X$. Then $D = \Sigma$ and $\phi(\Pi, P) \geq 2$ for each hyperplane $\Pi \subset \mathbf{P}^4$ and each point $P$ on $\Pi \cap X$. Now apply Bertini's theorem and choose a hyperplane $\Pi \subset X$ for which $\Pi \cap X$ is integral. Then, as $\phi(\Pi, P) \geq 2$ for each $P \in \Pi \cap X$, we conclude that each point on $\Pi \cap X$ lies on more than one line on $\Pi \cap X$ or on a line on $\Pi \cap X$ of multiplicity $\geq 2$. But then $\Pi \cap X$ would be of degree $\leq 2$ by lemma 7.2, thereby contradicting the assumption that $X$ is of degree $\geq 3$. We have thus that $R \neq X$. This completes the proof.

**Definition 7.6.** Let $k$ be an algebraically closed field and $X \subset \mathbf{P}^n$ be an integral hypersurface over $k$. Then $X \subset \mathbf{P}^n$ is said to be a scroll in $\mathbf{P}^{n-2}:s$, if there is a closed one-dimensional subscheme $C \subset F_{n-2}(X)$ such that the corresponding family $I' \subset X \times_k C$ of $(n-2)$-planes on $X$ is projected birationally onto $X$ (cf. e.g. [Rg$_2$]).

**Remarks 7.7.** (a) Let $k$ be an algebraically closed field of characteristic 0 and $X \subset \mathbf{P}^n$ be an integral hypersurface of dimension$\geq 2$ and degree$\geq 3$. Suppose that $X$ is covered by its linear divisors. Then it is easy to deduce from lemma 7.2 and the theorem of Bertini that there is only one $(n-2)$-plane on $X$ through a general point on $X \subset \mathbf{P}^n$. Hence $X \subset \mathbf{P}^n$ is a scroll in $\mathbf{P}^{n-2}:s$.

(b) It is possible to deduce theorem 7.5 from the results of B.Segre [Sg]. If e.g. $X \subset \mathbf{P}^4$ is an integral three-dimensional hypersurface of degree$\geq 3$, then it is proved in (op.cit.) that $X$ is a scroll in planes when dim $F_1(X) \geq 3$.

In the next result, we denote by $\Pi_a \subset \mathbf{P}^n$, the hyperplane given by $a_0 x_0 + \ldots + a_n x_n = 0$ for $a = (a_0, \ldots, a_n) \in \mathbf{P}^{n\vee}$.

**Corollary 7.8.** *Let $n \geq 4$ and $X \subset \mathbf{P}^n$ be a geometrically integral hypersurface over $\mathbf{Q}$ of degree $d \geq 3$ which is not covered by its linear divisors. Then there exists a form $\Phi(a_0, \ldots, a_n)$ over $\mathbf{Q}$ of degree $\ll_{d,n} 1$ such that the hyperplane section $\Pi_a \cap X$ is geometrically integral and not covered by its linear divisors when $\Phi(a_0, \ldots, a_n) \neq 0$.*

*Proof.* It is well known (cf. e.g. [BrSa], lemma 2.2.1) that there exists a form $\Phi_1(a_0, \ldots, a_n)$ over $\mathbf{Q}$ of degree $\ll_{d,n} 1$ such that $\Pi_a \cap X$ is geometrically integral when $\Phi_1(a_0, \ldots, a_n) \neq 0$. It is therefore sufficient to prove that there exists a form $\Phi_2(a_0, \ldots, a_n)$ over $\mathbf{Q}$ of degree $\ll_{d,n} 1$ such that $\Pi_a \cap X$ is not covered by its linear divisors when $\Phi_2(a_0, \ldots, a_n) \neq 0$.
Let $G$ (resp. $\mathbf{P}^{n\vee}$) be Grassmannian of linear $(n-3)$-subspaces $\Lambda$ of $\mathbf{P}^n$ (resp. hyperplanes $\Pi$ of $\mathbf{P}^n$) and $H$ be the Hilbert scheme of all hypersurfaces $X \subset \mathbf{P}^n$ of degree $d$. Let $I \subset G \times \mathbf{P}^{n\vee} \times H$ be the closed subscheme representing triples $(\Lambda, \Pi, X)$ such that $\Lambda \subset \Pi$ and $\Lambda \subset X$ and let $p: I \to \mathbf{P}^{n\vee} \times H$ be the projection map. Let $h_X$ be the Hilbert point of $X$ and $J \subset \mathbf{P}^{n\vee} \times H$ be the subset of points $(a, h_X)$ parametrising pairs $(\Pi_a, X)$ such that dim $p^{-1}(a, h_X) \geq 1$. Then $J$ is closed by upper semicontinuity of fibre dimension [G$_2$], Corollaire 13.1.5. If we introduce bihomogeneous coordinates for $\mathbf{P}^{n\vee} \times H$, then $J$ is defined by the vanishing of finitely many bihomogeneous polynomials $P_\lambda$. Now let $q: \mathbf{P}^{n\vee} \times H \to H$ be the projection. Then the fibre of $q_{|J}: J \to H$ at $h_X$ is a proper subset of $\mathbf{P}^{n\vee}$ by (7.5). Therefore, one of the bihomogeneous polynomials $P_\lambda$ must specialise to a non-zero form $\Phi_2(a_0, \ldots, a_n)$ of degree $\ll_{d,n} 1$ at $h_X$. Moreover, if $\Phi_2(a_0, \ldots, a_n) \neq 0$, then dim $F_{n-3}(\Pi_a \cap X) \leq 0$ as $F_{n-3}(\Pi_a \cap X) = p^{-1}(a, h_X)$. Hence $\Pi_a \cap X$ is not covered by its linear divisors, thereby completing the proof.

### 8. Integral and rational points on hypersurfaces

We are now in a position to generalise theorem 6.4 to affine hypersurfaces of dimension$\geq 2$.

**Theorem 8.1.** *Let $n \geq 3$ and $f(y_1, \ldots, y_n) \in \mathbf{Z}[y_1, \ldots, y_n]$ be a polynomial of degree $d \geq 4$ which is absolutely irreducible over $\mathbf{Q}$. Suppose that there are only finitely many affine linear $(n-2)$-subspaces over $\overline{\mathbf{Q}}$ on the variety $Y \subset \mathbf{A}^n$ defined by $f(y) = 0$. Then,*

$$n(Y\,;B) = O_{d,n,\varepsilon}(B^{n-2+\varepsilon})\,.$$

*Proof.* We shall use induction by dimension. The case $n=3$ is theorem 6.4. Suppose therefore that $n\geq 4$ and that 8.1 is true in lower dimensions. Let $F(x_0,\ldots,x_n) = x_0^d f(x_1/x_0,\ldots,x_n/x_0)$ and $X=\mathrm{Proj}(\mathbf{Q}[x_0,\ldots,x_n]/(F))$. Then $X$ is a geometrically integral hypersurface in $\mathbf{P}_\mathbf{Q}^n$. Also, $X$ cannot be covered by its linear divisors as there are only finitely many projective linear $(n-2)$-subspaces over $\overline{\mathbf{Q}}$ on $X$. Let $\Phi(a_0,\ldots,a_n)$ be a form in homogeneous coordinates $(a_0,\ldots,a_n)$ of $\mathbf{P}^{n\vee}$ as in corollary 7.8. There are $O_{d,n}(B^n)$ rational points $\boldsymbol{a}$ of $\mathbf{P}^{n\vee}$ of height $\leq B$ with $a_n\Phi(a_0,\ldots,a_n)=0$. We may thus find a primitive $(n+1)$-tuple $\boldsymbol{a}=(a_0,\ldots,a_n)$ of integers of height $<<_{d,n}1$ with $a_n\Phi(a_0,\ldots,a_n)\neq 0$ which will be fixed from now on. We shall also let $\Pi_b\subset\mathbf{P}^n$ be the hyperplane defined by the equation $bx_0+a_1x_1+\ldots+a_{n-1}x_{n-1}+a_nx_n=0$ and $\varphi(b)=\Phi(b,a_1,\ldots,a_n)$. Then the scheme-theoretic intersection $X_b=\Pi_b\cap X$ is not covered by its linear divisors when $\varphi(b)\neq 0$. There are thus for such $b$ only finitely many projective linear $(n-3)$-subspaces over $\overline{\mathbf{Q}}$ on $X_b$ by proposition 7.4.

Let $f_b(y_1,\ldots,y_{n-1}) = a_n^d f(y_1,\ldots,y_{n-1},-(b+a_1y_1+\ldots+a_{n-1}y_{n-1})/a_n)\in\mathbf{Z}[y_1,\ldots,y_{n-1}]$ and $Y_b\subset\mathbf{A}^{n-1}$ be the hypersurface defined by $f_b=0$. Then, $n(Y_b\,;B)$ can also be interpreted as the the number of $n$-tuples $\boldsymbol{y}=(y_1\ldots,y_n)$ in integers such that $y_1,\ldots,y_n\in[-B,B]$ and $f(y_1,\ldots,y_n) = b+a_1y_1+\ldots+a_{n-1}y_{n-1}+a_ny_n=0$. Therefore, $n(Y,B)\leq\Sigma_b\,n(Y_b,B)$ where $b$ runs over all integers in $[-cB,cB]$ for $c=|a_1|+\ldots+|a_n|<<_{d,n}1$.

If $\varphi(b)\neq 0$, then there are only finitely many affine linear $(n-3)$-subspaces over $\overline{\mathbf{Q}}$ on $Y_b$. Hence $n(Y_b\,;B) = O_{d,n,\varepsilon}(B^{n-3+\varepsilon})$ by the induction hypothesis. The total contribution to $n(Y;B)$ from all $b$ in $[-cB,cB]$ with $\varphi(b)\neq 0$ is therefore $O_{d,n,\varepsilon}(B^{n-2+\varepsilon})$.

To estimate the contribution from the hyperplane sections $Y_b\subset\mathbf{A}^{n-1}$ with $\varphi(b)=0$, we note that $\varphi$ cannot vanish identically since $\varphi(a_0)\neq 0$. There are thus at most $\deg\varphi<<_{d,n}1$ such hyperplane sections. For these we make use of the estimate $n(Y_b\,;B) = O_{d,n}(B^{n-2+\varepsilon})$ which follows from a sharper estimate in [P$_1$]. We thus get $O_{d,n}(B^{n-2+\varepsilon})$ in total for these $b$. This completes the proof.

The following result is an immediate corollary of theorem 8.1.

**Theorem 8.2.** *Let $Z\subset\mathbf{P}^{r+1}$ be a geometrically integral projective hypersurface over $\mathbf{Q}$ of degree $d\geq 4$. Suppose there are only finitely many $(r-1)$-planes over $\overline{\mathbf{Q}}$ on $Z$. Then,*

$$N(Z\,;B) = O_{d,r,\varepsilon}(B^{r+\varepsilon})\,.$$

*Proof.* Let $(z_1,\ldots,z_{r+2})$ be homogeneous coordinates for $\mathbf{P}^{r+1}$. Then $Z\subset\mathbf{P}^{r+1}$ is defined by some primitive integral form $f(z_1,\ldots,z_{r+2})$ of degree $d\geq 4$. Let $Y\subset\mathbf{A}^{r+2}$ be the affine cone of $Z\subset\mathbf{P}^{r+1}$ with affine coordinates $(y_1,\ldots,y_{r+2})$. Then $N(Z\,;B)\leq 2n(Y;B)$. As there are only finitely many linear $r$-subspaces over $\overline{\mathbf{Q}}$ on $Y$, we therefore obtain the estimate from theorem 8.1 for $n=r+2$.

The following proposition is sufficient for deducing corollary 0.2 from theorem 8.2.

**Proposition 8.3.** *Let $Z \subset \mathbf{P}^{r+1}$, $r \geq 2$ be an integral projective hypersurface of degree $\geq 3$ over an algebraically closed field of characteristic 0. Suppose that the singular locus of Z is of codimension at least three in Z. Then there are only finitely many $(r-1)$-planes on Z.*

*Proof.* If $r=2$, then 8.3 reduces to the well-known fact (cf. e.g. [Co]) that there are only finitely many lines on a non-singular surface in $\mathbf{P}^3$ of degree $\geq 3$. If $r \geq 3$, then we obtain from the theorem of Bertini [H$_2$], th.17.16 that there exists a 3-plane $\Gamma \subset \mathbf{P}^{r+1}$ such that $\Gamma \cap Z$ is non-singular. Hence $\Gamma \cap Z$ is not covered by its linear divisors as there are only finitely many lines on $\Gamma \cap Z$. But then $Z$ cannot be covered by its linear divisors since otherwise the intersections of these divisors with $\Gamma$ would cover $\Gamma \cap Z$. We therefore obtain the desired assertion from proposition 7.4.

### 9. Birational projections and rational points on projective varieties

We shall in this section generalise theorem 8.2 to projective varieties of arbitrary codimension. We shall thereby use the following result from [BwHeSa].

**Lemma 9.1.** *Let $X \subset \mathbf{P}^n$ be a geometrically integral projective variety over $\mathbf{Q}$ of dimension $r$ and degree $d \geq 2$. Then there exists a rational $(n-r-1)$-plane $\Lambda \subset \mathbf{P}^n$ disjoint from X and a linear morphism $\lambda: \mathbf{P}^n \setminus \Lambda \to \mathbf{P}^{r+1}$ over $\mathbf{Q}$ with the following properties.*

*(i) $H(\lambda(P)) \ll_{d,N} H(P)$ for rational points P on $\mathbf{P}^n \setminus \Lambda$.*
*(ii) There are at most d points on X in each fibre of $\lambda$.*
*(iii) There exists a proper closed subscheme T of Z such that $X \setminus \lambda^{-1}(T)$ is isomorphic to $Z \setminus T$ under $\lambda$ and such that $T \subset \mathbf{P}^{r+1}$ is given by $\ll_{d,n} 1$ equations of degree $\ll_{d,n} 1$.*
*(iv) The scheme-theoretic image Z of X under $\lambda$ is a geometrically integral hypersurface of degree d in $\mathbf{P}^{r+1}$.*

*Proof.* See §3 in [BwHeSa].

We are now in a position to prove the main theorem 0.1 of this paper.

**Theorem 9.2.** *Let $X \subset \mathbf{P}^n$ be a geometrically integral projective variety over $\mathbf{Q}$ of dimension r and degree $d \geq 4$. Suppose there are only finitely many $(r-1)$-planes over on X. Then,*

$$N(X ; B) = O_{d,n,\varepsilon}(B^{r+\varepsilon}).$$

*Proof.* Let $\lambda: X \to \mathbf{P}^{r+1}$ and $Z = \lambda(X) \subset \mathbf{P}^{r+1}$ be as in 9.1. Then it follows from (*i*) and (*ii*) in 9.1 that $N(X ; B) \leq dN(Z ; cB)$ for some $c \ll_{d,n} 1$. By 9.1(*iii*) we see that a linear $(r-1)$-plane on Z which is not contained in T is the image of a linear $(r-1)$-plane on X. There are thus only finitely many $(r-1)$-planes on Z. Also, deg $Z$ = deg $X = d \geq 4$ by 9.1(*iv*) so that $N(Z ; cB)$ $\ll_{d,n,\varepsilon} (cB)^{r+\varepsilon}$ by theorem 8.2. This completes the proof.

To deduce corollary 0.3 from theorem 9.2, we use the following result.

**Proposition 9.3.** *Let $Z \subset \mathbf{P}^{r+2}$, $r \geq 2$ be an integral r-dimensional intersection of two hypersurfaces of degree $d_1$ and $d_2$ over an algebraically closed field of characteristic 0.*

*Suppose that $d= d_1d_2 \geq 4$ and that the singular locus of X is of codimension at least three in X. Then there are only finitely many* ($r$-1)-*planes on Z.*

*Proof.* If $r=2$ , then 9.3 reduces to the known fact (cf. e.g. [Br], lemma 11) that there are only finitely many lines on a non-singular complete intersection surface $X \subset \mathbf{P}^4$ of degree $d_1d_2 \geq 4$. To handle the case $r>2$, we use the same Bertini argument as in the proof of corollary 8.3.

----------------------------------------------

Department of Mathematics , Chalmers University of Technology , S-412 96 Göteborg , Sweden
e-mail : salberg@math.chalmers.se